\documentclass{amsart}
\usepackage[latin1]{inputenc}

\usepackage{amssymb,amsmath}
\usepackage{amsthm}
\usepackage{paralist}
\usepackage{url}

\newtheorem{theorem}{Theorem}
\newtheorem{proposition}[theorem]{Proposition}
\newtheorem{lemma}[theorem]{Lemma}
\newtheorem{corollary}[theorem]{Corollary}

\newtheorem{problem}{Problem}
\newtheorem{fact}[theorem]{Fact}
\theoremstyle{remark}
\newtheorem{example}{Example}
\newtheorem{claim}{Claim}

\newcommand{\card}[1]{\ensuremath{\lvert{#1}\rvert}} % cardinality or length
\DeclareMathOperator{\ess}{ess} % essential arity
\DeclareMathOperator{\gap}{gap} % arity gap
\begin{document}

\title[Irreducible Boolean functions]{Irreducible Boolean functions}

\footnotetext[1]{Institut Sup\'erieur des Technologies M\'edicales de Tunis, $^2$University of Luxembourg, $^3$Universit\'e Claude-Bernard Lyon1}

\author[Moncef Bouaziz]{Moncef Bouaziz$^1$}
\address{Institut Sup\'erieur des Technologies M\'edicales de Tunis,\\
 9 avenue Dr Zouhair Essafi,\\
 1006 TUNIS, Tunisie  \\}
\email{Moncef.Bouaziz@istmt.rnu.tn}

%\footnotetext[2]{University of Luxembourg}
\author[Miguel Couceiro]{Miguel Couceiro$^2$}
\address{Department of Mathematics \\
University of Luxembourg \\
162a, avenue de la Fa\"iencerie \\
L-1511 Luxembourg}
\email{miguel.couceiro@uni.lu}

\author[Maurice Pouzet]{Maurice Pouzet$^3$}\thanks{The research of the first and last author has been supported by    CMCU Franco-Tunisien "Outils math\'ematiques pour l'informatique". }
\address {ICJ, Department of Mathematics, Universit\'e Claude-Bernard Lyon1, 43 Bd 11 Novembre 1918, 68622 Villeurbanne Cedex,
 France}
\email{pouzet@univ-lyon1.fr }


\date{\today}

\begin{abstract}  This paper is a contribution to the study of a quasi-order on the set $\Omega$ of Boolean functions, the \emph{simple minor} quasi-order. We look at the join-irreducible members of the resulting poset $\tilde{\Omega}$. Using a two-way correspondence between Boolean functions and hypergraphs, join-irreducibility translates into a combinatorial property of hypergraphs. We observe that among Steiner systems, those which yield join-irreducible members of $\tilde{\Omega}$ are the $-2$-monomorphic Steiner systems. We also describe the graphs which correspond to join-irreducible members of $\tilde{\Omega}$.
\end{abstract}
\maketitle

\section{Introduction}
Two approaches to define properties of Boolean functions have been considered in recent years; one in terms of
 functional equations \cite{EFHH}, and one other in terms of relational constraints \cite{Pippenger}. As it turned out,
these two approaches have the same expressive power in the sense that they specify exactly the same classes (or properties) of Boolean functions. The characterization of these classes was first obtained by Ekin, Foldes, Hammer and Hellerstein \cite{EFHH} who showed that equational classes of  Boolean functions can be completely described in terms of a quasi-ordering $\leq$ of the set $\Omega$ of all Boolean functions, called \emph{identification minor} in \cite{EFHH, H}, \emph{simple minor} in \cite {Pippenger, CP, CL1, CL2}, \emph{subfunction} in \cite {Zhegalkin}, and \emph{simple variable substitution} in \cite{CF1}.
 This quasi-order can be described as follows: for $f,g\in \Omega$, $g\leq  f $ if $g$ can be obtained from $f$ by identification of variables, permutation of variables or addition of dummy variables. As shown in \cite{EFHH}, equational classes of Boolean functions
 coincide exactly with the initial segments ${\downarrow }K=\{g\in \Omega: g\leq f, \textrm{ for some }f\in K\}$ of this quasi-order, or equivalently, they to correspond to antichains $A$ of Boolean functions in the sense that they constitute sets of the form $\Omega\setminus {\uparrow }A$. Similarly, those equational classes definable by finitely many equations where shown to correspond to finite antichains of Boolean functions. Since then, several investigations have appeared in this direction, to mention a few, see \cite{CF1, CF2, FP, Pippenger, Po}. 
 
 This correspondence to function class definability led to a greater emphasis on this quasi-ordering $\leq$ \cite{CP,CL1,CL2,CL3}.
As any quasi-order, the simple minor relation $\leq$ induces a partial order $\sqsubseteq$ on the set $\tilde{\Omega}$ made of equivalence classes of Boolean functions. Several properties of the resulting poset $(\tilde{\Omega},\sqsubseteq)$ were established in \cite{CP} where a classification of this poset was given in terms of equimorphism (two posets are \emph{equimorphic} if they are equivalent w.r.t. embeddings). Essentially, it was shown that $(\tilde{\Omega}, \sqsubseteq)$ has a sort of universal property among countable posets. 

 In this paper we are interested in the join-irreducible members of the poset $(\tilde{\Omega},\sqsubseteq)$, that is, those equivalence classes having a unique lower cover in $(\tilde{\Omega},\sqsubseteq)$. More precisely, we consider the problem of determining the join-irreducibles of this poset. Rather than taking a direct approach by looking into $(\tilde{\Omega},\sqsubseteq)$, we attack this problem by looking at hypergraphs. Indeed,  the fact that every Boolean function can be represented by a unique multilinear polynomial over the two-element field $ GF(2)$ allows to establish  a complete correspondence between Boolean functions and hypergraphs. This correspondence is given as follows.
 For any hypergraph $\mathcal{H} = (V, \mathcal{E})$ we associate the polynomial $P_{\mathcal{H}}\in GF(2)[x_i, i\in V]$ given by $P_{\mathcal{H}}=\underset{E\in \mathcal{E}}{\sum} \underset{i\in E}{\prod}x_i$. Conversely, every polymomial $P\in GF(2)[x_i, i\in V]$ is of the form $P=P_{\mathcal{H}}$ where $\mathcal{H} = (V, \mathcal{E})$ and $\mathcal{E}$ is the set of hyperedges corresponding to the monomials of $P$.
%  Thus the mapping $\mathcal{H} \mapsto P_{\mathcal{H}}$ establishes a two-way correspondence between Boolean functions and hypergraphs.  

To work in complete analogy with the Boolean function setting, we still need to mimic the simple minor relation in the realm of hypergraphs. This is acheived through the notion of quotient map. Say that a map $h'\colon V'\to V$ is a \emph{quotient map} from $\mathcal{H}' = (V', \mathcal{E}')$ to $\mathcal{H} = (V, \mathcal{E})$ if for every $E\subseteq V$, $E\in \mathcal{E}$ if and only if $\card{\{E'\in \mathcal{E}':{h'}(E')=E}$ is odd.
For two hypergraphs $\mathcal{H}'$ and $\mathcal{H}$, set $\mathcal{H}\preceq \mathcal{H}'$ if there is a quotient map from $\mathcal{H}' $ to $\mathcal{H}$. As we are going to see $\preceq$ constitutes a quasi-order between hypergraphs and two hypergraphs are related by $\preceq$ if and only if the corresponding Boolean functions are related by $\leq$ (see Lemma \ref{quasi} and Theorem \ref{correspondence}, resp.).  
The fact that a Boolean function (or more accurately an equivalence class) is join-irreducible translates to a combinatorial property of the corresponding hypergraph. A description of  all hypergraphs satisfying this property eludes us.  But, as we will observe, among these hypergraphs some have been intensively studied for other purposes. The basic examples are the non-trivial hypergraphs whose automorphism group is $2$-set transitive. We show that Steiner systems which yield join-irreducible members of the poset $(\tilde{\Omega},\sqsubseteq)$ are exactly those which are $-2$-\emph{monomorphic} in the sense that the induced hypergraphs obtained by deleting any pair of two distinct vertices are isomorphic (Theorem \ref {thm:steiner}). Among Steiner triple systems those with a flag-transitive automorphism group enjoy this property. We do not know if there are other.
We also describe those graphs corresponding to join-irreducible members of $(\tilde{\Omega},\sqsubseteq)$ (Theorem \ref{thm:jigraphs}). %\section{Posets}

%Join-irreducibles, covers, etc.

\section{Boolean functions}

A \emph{Boolean function} is simply a mapping $f\colon \{0,1\}^n\to \{0,1\}$ where $n\geq 1$ and called the \emph{arity} of $f$. The simplest examples of $n$-ary Boolean functions are the \emph{projections}, i.e., mappings $(a_1,\ldots ,a_n)\mapsto a_i$, for $1\leq i\leq n$ and $a_1,\ldots ,a_m\in \{0,1\}$, and which we also refer to as \emph{variables}. 
For each $n\geq 1$, we denote by $\Omega^{(n)}=\{0,1\}^{\{0,1\}^n}$ the set of all $n$-ary Boolean functions and we denote by $\Omega=\bigcup_{n\geq 1}\Omega^{(n)}$ the set of all Boolean functions. 

A variable $x_i$ is  an \emph{essential} variable of $f$ if  $f$ \emph{depends} on its $i$-th argument, that is if  there are $a_1,\ldots ,a_{i-1},a_{i+1},\ldots ,a_{n}\in\{0,1\}$ such that
 the unary function $f(a_1,\ldots ,a_{i-1},x_i,a_{i+1},\ldots ,a_{n})$ is nonconstant. By \emph{essential arity} of a function $f\in\Omega^{(n)}$, denoted $\ess f$, we simply mean the number of its essential variables. 
 For instance, constant functions are exactly those functions with essential arity $0$. Functions with essential arity $1$ are either projections or negated projections.

\subsection{Simple minors and irreducible Boolean functions}

A Boolean function $g \in\Omega^{(m)}$ is said to be a \emph{simple minor} of a Boolean function $f\in \Omega^{(n)}$
if there is a mapping $\sigma \colon \{1, \ldots, n\} \to \{1, \ldots, m\}$ such that
$$g = f(x_{\sigma(1)}, \ldots, x_{\sigma(n)}).$$
If $\sigma$ is not injective, then we speak of \emph{identification of variables.} If $\sigma$ is not surjective, then we speak of \emph{addition of inessential variables.} If $\sigma$ is a bijection, then we speak of \emph{permutation of variables.} In fact, these Mal\'cev operations are sufficient to completely describe the simple minor relation.

%\begin{fact} A Boolean function $g \in\Omega^{(m)}$ is a simple minor of $f\in \Omega^{(n)}$ if and only if $g \circ I_c \subseteq f \circ I_c$, where
%\[
%f \circ I_c = \{f(p_1, \ldots, p_n) : \text{for projections $p_1, \ldots, p_n$ of the same arity}\}.
%\]
%\end{fact}

\begin{fact} 
 The simple minor relation between Boolean functions is a quasi-order. 
\end{fact}

Let $\leq$ denote the simple minor relation on the set $\Omega $ of all Boolean functions. If $g \leq f$ and $f \leq g$, then we say that $f$ and $g$ are \emph{equivalent,} denoted $f \equiv g$. The equivalence class of $f$ is denoted by $\tilde{f}$. If $g \leq f$ but $f \not\leq g$, then we use the notation $g < f$. The arity gap of $f$, denoted $\gap f$, is defined by $\gap f=min\{\ess f - \ess g: g<f\}$. Note that equivalent functions may differ in arity, but not in essential arity nor in arity gap.

\begin{fact}\label{essfact} If $g \leq f$, then $\ess g \leq \ess f$, with equality if and only if $g \equiv f$.
\end{fact}

Let $(\tilde{\Omega}, \sqsubseteq)$ denote the poset made of equivalence classes of Boolean functions associated with the simple minor relation, that is, $\tilde{\Omega}=\Omega /\equiv$ together with the partial order $\sqsubseteq$ given by $\tilde{g}\sqsubseteq \tilde{f}$ if and only if $g\leq f$. Several properties of this poset were established in \cite{CP}. For example, Fact \ref{essfact} implies that each principal initial segment ${\downarrow \tilde{f}}=\{\tilde{g}:\tilde{g}\sqsubseteq \tilde{f}\}$ is finite. This means that $(\tilde{\Omega}, \sqsubseteq)$ decomposes into levels $\tilde{\Omega}_{0}, \dots  ,\tilde{\Omega}_{n}, \dots$, where $\tilde{\Omega}_{n}$ is the set of minimal elements of $\tilde{\Omega}\setminus \bigcup \{ \tilde{\Omega}_{m}: m<n\}$. 

\begin{fact} If $\ess g=n$ and $\tilde{g}$ is covered by $\tilde{f}$, then $\ess f= n + \gap f$.
\end{fact}

This fact and Salomaa's result \cite{Salomaa} which asserts that the arity gap of Boolean functions is at most $2$, imply that each level of $(\tilde{\Omega}, \sqsubseteq)$ is finite.
 For instance, the first level $\tilde{\Omega}_{0}$ comprises four equivalence classes, namely, those of constant $0$ and $1$ functions, and those of projections and negated projections. These four classes induce a partition of $(\tilde{\Omega}, \sqsubseteq)$ into four different blocks with no comparabilities in between them. For further background see \cite{CP}.

We say that Boolean function $f$ is \emph{irreducible} if there is $f'\in \Omega$ such that 
\begin{enumerate}[(i)]
\item $f'<f$, and 
\item for every $g \in \Omega$,  if $g<f$, then $g \leq f'$. 
\end{enumerate}

\begin{fact}
 A Boolean function $f$ is irreducible if and only if $\tilde{f}$ has a unique lower cover, i.e., $\tilde{f}$ is join-irreducible in $(\tilde{\Omega}, \sqsubseteq)$.
\end{fact}

To illustrate, consider the binary conjunction $x_1\wedge x_2$, the binary disjunction $x_1\vee x_2$ and the composite $(x_1\vee x_2) \wedge (x_3\vee x_4)$. Both the binary conjuction and disjunction constitute irreducible functions since they have, up to equivalence, a unique strict minor, namely, a projection. This uniqueness clearly extends to any conjunction and disjunction of $n\geq 2$ variables, showing that any of the latter functions also constitute irreducible functions. But this is not the case for the composite $(x_1\vee x_2) \wedge (x_3\vee x_4)$. Indeed, $x_1\wedge x_3 \vee x_1\wedge x_4, x_1\wedge x_4 \vee x_2 \leq (x_1\vee x_2) \wedge (x_3\vee x_4)$, but $x_1\wedge x_3 \vee x_1\wedge x_4 \not \equiv x_1\wedge x_4 \vee x_2$. These observations lead to the following problem.
  
\begin{problem}\label{mainproblem} Describe the irreducible Boolean functions.
 \end{problem}

%\begin{fact} If $\ess g=n$ and $\tilde{g}$ is covered by $\tilde{f}$, then $\ess f= n + \gap f$.
%\end{fact}

%Using this fact and Salomaa's \cite{Salomaa} result which asserts that the arity gap of Boolean functions is at most $2$,  we were able to show that each level of $(\tilde{\Omega}, \sqsubseteq)$ is finite. For example, the first level $\tilde{\Omega}_{0}$ comprises four equivalence classes, namely, those of constant $0$ and $1$ functions, and those of projections and negated projections. In fact, these four classes induce a partition of $(\tilde{\Omega}, \sqsubseteq)$ into four different blocks with no comparabilities in between them. (See \cite{CP} for further background and see \cite{CL1} for a complete classification of Boolean functions according to their arity gap. See also \cite{CL2} for further extensions concerning the arity gap of finite functions.)

\subsection{Boolean functions as polynomials}

In this subsection, we view $\{0,1\}$ as endowed with the two-element field structure, $\{0, 1\} = GF(2)$,
as well as with the lattice structure where $0 < 1$. Consider the commutative ring $GF(2)[x_1,\ldots ,x_n]$ of multilinear polynomials in $n$ indeterminates. Each of these polynomials is characterized by the fact that each monomial is a product of distinct indeterminates.

% In the sequel, we shall regard monomials has maps from ${\bf n}=\{1,\ldots , n\}$ to $\{0,1\}$ where each such map $h$ is denoted by $\underset{i\in {\bf n}}{\prod }x_i^{h(i)}$. For instance, the monomial $1$ is the constant map $0$.

To each polynomial $P\in GF(2)[x_1,\ldots ,x_n]$ corresponds an $n$-ary Boolean function $f_P \colon \{0,1\}^n\to \{0,1\}$ which is given as the evaluation of $P$, that is, for every
$(a_1,\ldots ,a_n)\in \{0,1\}^n$, $f_P(a_1,\ldots ,a_n)=P(a_1,\ldots ,a_n)$. The function $f_P$ is said to be \emph{represented} by $P$, and $P$ is said to be 
\emph{Zhegalkin} (or \emph{Reed--Muller}) polynomial of $f_P$ \cite{Muller,Reed,Zhegalkin}. As it is well-known every Boolean function can be represented in this way.

\begin{theorem} Every Boolean function $f\colon \{0,1\}^n\to \{0,1\}$, $n\geq 1$, is uniquely represented by a multilinear polynomial $P\in GF(2)[x_1,\ldots ,x_n]$.
\end{theorem}

This result allows to work with polynomials rather than Boolean functions. This approach turns out to be quite useful when studying the poset  $(\tilde{\Omega}, \sqsubseteq)$. For instance, as we mentioned the four equivalence classes in $\tilde{\Omega}_{0}$, namely, those represented by $0,1,x_1$ and $x_1+1$ induce a partition of $\tilde{\Omega}$ into different blocks with no comparabilities in between them. As it is easy to verify, above the equivalence classes represented by the constant polynomials $0$ or $1$ we have the equivalence classes of those functions whose Zhegalkin polynomials are the sum of an even number of nonconstant monomials plus $0$ or $1$, respectively, and above the equivalence classes represented by $x_1$ or $x_1+1$ we have the equivalence classes of those functions whose Zhegalkin polynomials are the sum of an odd number of nonconstant monomials plus $0$ or $1$, respectively.

\begin{corollary}
 A variable $x_i$ is essential in $f\in \Omega^{(n)}$ if and only if $x_i$ appears in the Zhegalkin polynomial of $f$. In particular, $\ess f$ is the number of variables appearing in the Zhegalkin polynomial of $f$.
\end{corollary}

Thus, in the case of polynomial expressions, to describe the simple minor relation we only need to consider identification and permutation of essential variables, since the operation of addition of inessential variables produces the same polynomial representations. Moreover, from Fact \ref{essfact} it follows that the strict minors of a given function $f$ have Zhegalkin polynomials with strictly less variables, and that the Zhegalkin polynomials of functions equivalent to $f$ are obtained from the Zhegalkin polynomial of $f$ by permutation of its variables. 
For further developements see \cite{CL3}.

We finish this section with a complete classification of Boolean functions according to their arity gap, and which we shall make use of in the following sections.

\begin{theorem}\label{booleangap}\emph{(In \cite{CL1}:)}
Let $f \colon \{0, 1\}^n \to \{0, 1\}$ be a Boolean function with at least two essential variables. Then the arity gap of $f$ is two if and only if its Zhegalkin polynomial is equivalent to one of the following:
\begin{compactenum}
\item $x_1 + x_2 + \dots + x_m + c$ for some $m \geq 2$,
\item $x_1 x_2 + x_1 + c$,
\item $x_1 x_2 + x_1 x_3 + x_2 x_3 + c$,
\item $x_1 x_2 + x_1 x_3 + x_2 x_3 + x_1 + x_2 + c$,
\end{compactenum}
where $c \in \{0,1\}$. Otherwise the arity gap of $f$ is one.
\end{theorem}

\section{Boolean functions and hypergraphs}

By an \emph{hypergraph} we simply mean a pair $\mathcal{H} = (V, \mathcal{E})$ where $V$ is a finite nonempty set whose elements are called \emph{vertices}, and where $\mathcal{E}$ is a collection of subsets of $V$ called \emph{hyperedges}. If $V$ has $n$ elements, then we set $V=\{1,\ldots ,n\}$ and we write $[V]^m$ to denote the set of $m$-element subsets of $V$.

Let $\mathcal{H} = (V, \mathcal{E})$ be an hypergraph with $n$ vertices. To such an hypergraph $\mathcal{H}$ we associate a polynomial $P_{\mathcal{H}}\in GF(2)[x_i, i\in V]$ which is given by $P_{\mathcal{H}}=\underset{E\in \mathcal{E}}{\sum} \underset{i\in E}{\prod}x_i$.

\begin{example}
Let $\mathcal{H}_1 = (\{1,2,3\}, \emptyset)$,   $\mathcal{H}_2 = (\{1,2,3\},\{\{1,2\}, \emptyset\})$ and $\mathcal{H}_3 = (\{1,2,3\}, \{\{1,2\}, \{1,3\}, \{2,3\}\} )$. Then $\mathcal{P}_{\mathcal{H}_1}= 0$, $\mathcal{P}_{\mathcal{H}_2}= x_1x_2+1$ and $\mathcal{P}_{\mathcal{H}_3}= x_1x_2+x_1x_2+x_2x_3$, respectively.
\end{example}

Conversely, it is clear that to each polymomial $P\in GF(2)[x_1,\ldots ,x_n]$ is associated an hypergraph 
$\mathcal{H}_P = (V, \mathcal{E})$ where $V=\{1,\ldots ,n\}$ and $\mathcal{E}$ is the set of hyperedges corresponding to the monomials of $P$. Thus, using the two-way correspondence between Boolean functions and polynomials over $GF(2)$, we have the following.

\begin{theorem} For each Boolean function $f\colon \{0,1\}^n\to \{0,1\}$, $n\geq 1$, there is a unique hypergraph $\mathcal{H} = (V, \mathcal{E})$, $V=\{1,\ldots ,n\}$, such that $f=f_{P_{\mathcal{H}}}$.
\end{theorem}

For the sake of simplicity, let $f_{\mathcal{H}}$ denote the function $f_{P_{\mathcal{H}}}$ determined by $\mathcal{H}$.

\subsection{Simple minors of hypergraphs}

Let $\mathcal{H} = (V, \mathcal{E})$ and $\mathcal{H}' = (V', \mathcal{E}')$ be two hypergraphs and let $h'\colon V'\to V$ be a map. For each $E\subseteq V$, set ${h'}^{-1}[E]=\{E'\in \mathcal{E}': h'(E')=E\}$, where $h'(E')=\{h'(i'): i'\in E'\}$.
The map $h'$ is said to be a \emph{quotient map from $\mathcal{H}'$ to $\mathcal{H}$}, denoted $h'\colon \mathcal{H}'\to \mathcal{H}$, if for every 
$E\subseteq V$, the following condition holds: $E\in \mathcal{E}$ if and only if the cardinality $\card{{h'}^{-1}[E]}$ is odd.
We say that an hypergraph $\mathcal{H}$ is a \emph{simple minor} of an hypergraph $\mathcal{H}'$, denoted $\mathcal{H}\preceq \mathcal{H}'$, if there is a quotient map from $\mathcal{H}'$ to $\mathcal{H}$.

%To illustrate, let $\mathcal{H} = (V, \mathcal{E})$ be an hypergraph with $V=\{1,\ldots ,n\}$. For $i,j\in V$, consider the hypergraph $\mathcal{H}_{i=j} = (V_{i=j}, \mathcal{E}_{i=j})$ given as follows:
%\begin{itemize}
% \item $V_{i=j}=V\setminus \{i\}$, and 
%\item for each $E \subseteq V_{i=j}$, $E \in \mathcal{E}_{i=j}$ if either
%\begin{enumerate}[(i)]
% \item $E \in \mathcal{E}$ and $\{i,j\}\cap E=\emptyset$, or
%\item \label{oddcondition} $j\in E$ and among $E$, $(E\setminus \{j\})\cup \{i\}$ and $E\cup \{i\}$, either one or the three sets belong to $\mathcal{E}$.
%\end{enumerate}
%\end{itemize}
%Note that the condition (\ref{oddcondition}) guarantees that the map $h\colon V\to V_{i=j}$ defined by $h(i)=h(j)=j$ and $h(k)=k$, for each $k\neq i,j$, constitutes a quotient map from $\mathcal{H}$ to $\mathcal{H}_{i=j}$, thus showing that $\mathcal{H}_{i=j}$ is a simple minor of $\mathcal{H}$.
%Note also that $\mathcal{H}_{i=j}\cong \mathcal{H}_{j=i}$.

To illustrate, let $\mathcal{H} = (V, \mathcal{E})$ be an hypergraph with $V=\{1,\ldots ,n\}$. Let $e=\{i,j\}$, $i,j\in V$, and fix $l_e\not \in V$. Consider the hypergraph $\mathcal{H}_{e} = (V_{e}, \mathcal{E}_{e})$ given as follows: $V_{e}=(V\setminus e)\cup \{l_e\}$ and for each $E \subseteq V_{e}$, we have $E \in \mathcal{E}_{e}$ if either
\begin{enumerate}[(i)]
 \item $E \in \mathcal{E}$ and $e\cap E=\emptyset$, or
\item \label{oddcondition} $l_e\in E$ and among $(E\setminus \{l_e\})\cup e$, $(E\setminus \{l_e\})\cup \{i\}$ and $(E\setminus \{l_e\})\cup \{j\}$, either one or the three sets belong to $\mathcal{E}$.
\end{enumerate}
Note that the condition (\ref{oddcondition}) guarantees that the map $h\colon V\to V_{e}$ defined by $h(i)=h(j)=l_e$ and $h(k)=k$, for each $k\neq i,j$, constitutes a quotient map from $\mathcal{H}$ to $\mathcal{H}_{e}$, thus showing that $\mathcal{H}_{e}$ is a simple minor of $\mathcal{H}$.

 \begin{lemma}\label{quasi}
  The simple minor relation between hypergraphs is a quasi-order.
 \end{lemma}

\begin{proof}
 Let $\mathcal{H} = (V, \mathcal{E})$, $\mathcal{H}' = (V', \mathcal{E}')$ and $\mathcal{H}^{''} = (V^{''}, \mathcal{E}^{''})$ be hypergraphs such that $\mathcal{H}\preceq \mathcal{H}^{'}\preceq \mathcal{H}^{''}$. Let $h'\colon \mathcal{H}'\to \mathcal{H}$ and $h^{''}\colon \mathcal{H}^{''}\to \mathcal{H}'$ be the corresponding quotient maps. We claim that $h=h'\circ h^{''}$ is a quotient map from $\mathcal{H}^{''}$ to 
$\mathcal{H}^{'}$. 

Let $E\subseteq V$. The set $h^{-1}[E]=\{E^{''}\in \mathcal{E}^{''}:h^{''}(E)=E \}$ decomposes into two sets, namely, 
$A=\bigcup \{{h^{''}}^{-1}[E']:E'\in \mathcal{E}', h'(E')=E\}$ and 
$B=\bigcup \{{h^{''}}^{-1}[E']:E'\not\in \mathcal{E}', h'(E')=E\}$.
Now $A$ is a disjoint union of sets of odd size and $B$ is a disjoint union of sets of even size and hence, $B$ has even size. Thus the parity of $\card{h^{-1}[E]}$ is the same as the parity of $\card{A}$ which, in turn, is the same as the parity of $\card{{h'}^{-1}[E]}$. Since $E\in \mathcal{E}$ if and only if $\card{{h'}^{-1}[E]}$ is odd, the proof is now complete.
\end{proof}

The following theorem establishes the connection between the simple minor relation on Boolean functions and the simple minor relation on hypergraphs.

\begin{theorem}\label{correspondence} Let $\mathcal{H} = (V, \mathcal{E})$ and $\mathcal{H}' = (V', \mathcal{E}')$ be two hypergraphs, with  $V=\{1,\ldots ,n\}$ and 
$V'=\{1,\ldots ,m\}$, respectively. Then $\mathcal{H}\preceq \mathcal{H}'$ if and only if $f_{\mathcal{H}}\leq f_{\mathcal{H}'}$.
\end{theorem}

\begin{proof} Suppose first that $\mathcal{H}\preceq \mathcal{H}'$ and let $h'\colon \mathcal{H}'\to \mathcal{H}$ be the quotient map. Define 
$\sigma \colon \{1, \ldots, m\} \to \{1, \ldots, n\}$ by $\sigma (i')=h(i')$, for every $i'\in V'=\{1,\ldots ,m\}$. To verify $f_{\mathcal{H}}\leq f_{\mathcal{H}'}$, we only have to show that 
$P_{\mathcal{H}}= P_{\mathcal{H}'}(x_{\sigma(1)}, \ldots, x_{\sigma(m)})$. Indeed, 
\begin{multline*}
P_{\mathcal{H}'}(x_{\sigma(1)}, \ldots, x_{\sigma(m)})= \underset{E'\in \mathcal{E}'}{\sum} \underset{i'\in E'}{\prod}x_{\sigma(i')}=
\underset{E'\in \mathcal{E}'}{\sum} \underset{i'\in E'}{\prod}x_{h(i')}= 
\\
\underset{E'\in \mathcal{E}'}{\sum} \underset{i\in h'(E')}{\prod}x_{i}^{\card{{h'}^{-1}(i)}}.
\end{multline*}
 
 Since $x_j^2=x_j$, we have that $\underset{i\in h'(E')}{\prod}x_{i}^{\card{{h'}^{-1}(i)}}=\underset{i\in h'(E')}{\prod}x_{i}$. 
Hence, $$ P_{\mathcal{H}'}(x_{\sigma(1)}, \ldots, x_{\sigma(m)})=\underset{E\subseteq V}{\sum}\, \underset{h'(E')=E}{\sum}\, \underset{i\in h'(E')}{\prod}x_{i} .$$
Now each term $\underset{h'(E')=E}{\sum}\, \underset{i\in h'(E')}{\prod}x_{i}$ is different to $0$ if and only if $\card{{h'}^{-1}[E]}$ is odd, that is, $E\in \mathcal{E}$. In other words, $ P_{\mathcal{H}'}(x_{\sigma(1)}, \ldots, x_{\sigma(m)})=P_{\mathcal{H}}$.

Now suppose that $f_{\mathcal{H}}\leq f_{\mathcal{H}'}$ and let $\sigma \colon \{1, \ldots, m\} \to \{1, \ldots, n\}$ be a map such that 
 $P_{\mathcal{H}}= P_{\mathcal{H}'}(x_{\sigma(1)}, \ldots, x_{\sigma(m)})$. Define $h'\colon V'\to V$ by $h(i')=\sigma (i')$, for every $i'\in V'=\{1,\ldots ,m\}$. Let $E\subseteq V$. We have that $E\in \mathcal{E}$ if and only if $\underset{i\in E}{\prod}x_i$ is a monomial of $P_{\mathcal{H}}$. Since $P_{\mathcal{H}}= P_{\mathcal{H}'}(x_{\sigma(1)}, \ldots, x_{\sigma(m)})$, the latter condition holds if and only if the number of monomials $\underset{i'\in E'}{\prod}x_{i'}$ of $P_{\mathcal{H}'}$ such that $\sigma(E')=E$ is odd. In other words, $E\in \mathcal{E}$ if and only if $\card{{h'}^{-1}[E]}$ is odd. This shows that $h'$ constitutes a quotient map from $\mathcal{H}'$ to $\mathcal{H}$. 
\end{proof}

\subsection{Conditions for irreducibility}

Let $\mathcal{H} = (V, \mathcal{E})$ and $\mathcal{H}' = (V', \mathcal{E}')$ be two hypergraphs. A map  $\varphi \colon V\to V'$ is said to be an \emph{isomorphism} from $\mathcal{H}$ onto $\mathcal{H}'$ if $\varphi$ is bijective and for every $E\subseteq V$, $E\in \mathcal{E}$ if and only if $\varphi(E) \in \mathcal{E}'$. Two hypergraphs  $\mathcal{H} $ and $\mathcal{H}' $ are said to be \emph{isomorphic}, denoted $ \mathcal{H} \cong \mathcal{H}' $, if there is an isomorphism $\varphi $ from $\mathcal{H}$ onto $\mathcal{H}'$. 
If $\mathcal{H} =\mathcal{H}' $, then $\varphi$ is called an \emph{automorphism} of $\mathcal{H}$. The group made of automorphisms of $\mathcal{H}$ is denoted by $Aut(\mathcal{H})$.

Let $\mathcal{H} = (V, \mathcal{E})$ be an hypergraph and let $\overline{V}=\bigcup \mathcal{E}$. For $e,e'\in [\overline{V}]^2$, define $e\approx e'$ if $\mathcal{H}_e\cong \mathcal{H}_{e'}$. Obviously, $\approx$ is an equivalence relation.

\begin{lemma}\label{keylemma}
 Let $\mathcal{H} = (V, \mathcal{E})$ be an hypergraph. Then
$f_{\mathcal{H}}$ is irreducible if and only if $\card{\overline{V}}\geq 2$ and there is an equivalence class $\mathcal{C}$ of $\approx$ such that, for every $e'\in [\overline{V}]^2\setminus \mathcal{C}$ and $e\in \mathcal{C}$, $\ess f_{\mathcal{H}_{e'}}<\ess f_{\mathcal{H}_e}$.
\end{lemma}

In the search for hypergraphs $\mathcal{H} = (V, \mathcal{E})$ determining irreducible Boolean functions, Lemma \ref{keylemma} invites us to look at differences $\ess f_{\mathcal{G}}-\ess f_{\mathcal{G}_{e}}$, especially, when $\ess f_{\mathcal{G}}-\ess f_{\mathcal{G}_{e}}>1$. For the latter to occur, there are two possibilities:
\begin{enumerate}[(i)]
 \item the vertex $l_e$ becomes isolated and this is the case if and only if, for every $F$ disjoint from $e$, the number of $e'\subseteq V$ such that $\emptyset \neq e'\subseteq e$ and $e'\cup F \in \mathcal{E}$, is even, or
\item another vertex, say $i\in V$, becomes isolated and this is the case if and only if, for every $e'\in \mathcal{E}$, if $i\in e'$ then $e\cap e'\neq \emptyset$ and there is $e''\in \mathcal{E}$ such that $i\in e''$ and $e'\setminus e=e''\setminus e$.
\end{enumerate}  
Note that if  $\ess f_{\mathcal{G}}-\ess f_{\mathcal{G}_{e}}>1$, for every $e\in [{V}]^2$, then by Theorem \ref{booleangap} it follows that $\ess f_{\mathcal{G}}-\ess f_{\mathcal{G}_{e}}=2$, for every $e\in [{V}]^2$, and $\mathcal{G}$ determines a function which is equivalent to one listed in Theorem \ref{booleangap}.

As an immediate consequence of Lemma \ref{keylemma}, we get the following criterion for irreducibility.

\begin{corollary} Let $\mathcal{H} = (V, \mathcal{E})$ be an hypergraph. If $\card{\overline{V}}\geq 2$ and, for every $e,e'\in [\overline{V}]^2$, we have $e\approx e'$, then $f_{\mathcal{H}}$ is irreducible.
\end{corollary}

A group $G$ acting on a set $V$ is \emph{$2$-set transitive} if for every $e,e'\in [V]^2$, there is some $g\in G$ such that $g(e)=e'$.

%supprime \underline sur V$
\begin{corollary} Let $\mathcal{H} = (V, \mathcal{E})$ be an hypergraph. If $\card{V}\geq 2$, $\bigcup \mathcal{E}= V$ and $Aut(\mathcal{H})$ is $2$-set transitive, then $f_{\mathcal{H}}$ is irreducible.
\end{corollary}
\begin{proof}
 Let $\varphi \in Aut(\mathcal{H})$. Take $e\in [{V}]^2$ and let $e'=\varphi(e)\in [{V}]^2$. Consider the mapping $\overline{\varphi}:V_e\to V_{e'}$ defined by $\overline{\varphi}(l_e)=l_{e'}$ and $\overline{\varphi}(i)=i$ for every $i\in V_e\setminus \{l_e\}$. Clearly, $\overline{\varphi}$ constitutes the desired isomorphism from $\mathcal{H}_e$ to $\mathcal{H}_{e'}$.
\end{proof}

\begin{problem}
 For which hypergraphs $\mathcal{H}$
\begin{enumerate}[(i)]
 \item $Aut(\mathcal{H})$ is $2$-set transitive?
\item $\mathcal{H}_e\cong \mathcal{H}_{e'}$, for every $e,e'\in [{V}]^2$?
\end{enumerate}
\end{problem}

Examples are given  in the next subsection. 
\subsection{Steiner Systems}

Let $\mathcal{H}=(V, \mathcal{E})$ be an hypergraph. We say that the hypergrapph $\mathcal{H}$ is a \emph{$2-(n,k,\lambda)$ design} if $\card{V}=n$, $\mathcal{E}\subseteq [V]^k$, and for every $e\in [V]^2$, $\card{\{E\in \mathcal{E}: e\subseteq E\}}=\lambda$.
 If $\lambda =1$, then we say that $\mathcal{H}$ is a \emph{Steiner system} and, in addition, if $k=3$, then we say that $\mathcal{H}$ is a \emph{Steiner triple system}.

For each $e\in [V]^2$, set $\mathcal{H}_{\setminus e}=(V\setminus e, \mathcal{E}\cap [V\setminus e]^2)$.
If for every $e,e'\in [V]^2$, $\mathcal{H}_{\setminus e}\cong \mathcal{H}_{\setminus e'}$, then we say that 
$\mathcal{H}$ is \emph{$-2$-monomorphic}. The following theorem reveals a connection between the notion of $-2$-monomorphicity and irreducibility in the case of Steiner systems.

\begin{theorem} \label{thm:steiner}Let $\mathcal{H}=(V, \mathcal{E})$ be a Steiner system. The following are equivalent:
 \begin{enumerate}[(i)]
  \item $f_{\mathcal{H}}$ is irreducible;
\item $\mathcal{H}_{ e}\cong \mathcal{H}_{ e'}$, for every $e,e'\in [V]^2$;
\item $\mathcal{H}$ is $-2$-monomorphic.
 \end{enumerate}
\end{theorem}

\begin{problem}
For a Steiner triple system $\mathcal{H}=(V, \mathcal{E})$ does the following hold:
$\mathcal{H}$ is $-2$-monomorphic if and only if $Aut(\mathcal{H})$ is $2$-set transitive?
\end{problem}
Note that the automorphism group of a Steiner  systems is flag-transitive whenever it is $2$-set transitive. The converse holds for Steiner triple systems. There are several deep results about Steiner systems with a $2$-transitive or a flag transitive automorphism group (see the survey by Kantor \cite{kantor}). For example, any Steiner triple system with a $2$-transitive automorphism group must be a projective space over $GF(2)$ or an affine space over $GF(3)$\cite{hall, key}. The notion of monomorphy (with some of its variations) is due to R. Fra{\"\i}ss\'e. His book \cite{fraissetr} contains some important results concerning this notion. 

\section{Join-irreducible graphs}\label{irreduciblegraphs}

In this subsection, we give an answer to Problem \ref{mainproblem} in the particular case of functions which are determined by simple graphs, that is graphs $\mathcal{G} = (V, \mathcal{E})$ where $\mathcal{E}\subseteq [{V}]^2$. As in the case of  hypergraphs, if we remove the  isolated vertices of  $\mathcal G$, the resulting graph $\check{\mathcal G}$ yields an equivalent function. In the sequel, when we speak of a \emph{join-irreducible} graph we simply mean a graph $\mathcal{G}$ such that $\tilde{f_{\mathcal{G}}}$ is join-irreducible. Note that  $\mathcal G$ is join-irreducible if and only if $\check {\mathcal G}$ is join-irreducible; note also that a join-irreducible graph must have at least one edge.

Given a graph $\mathcal{G} = (V, \mathcal{E})$, we write $i\sim j$ if $\{i,j\}\in \mathcal{E}$. Set $V(i)=\{j\in V: i\sim j\}\cup \{i\}$. The \emph{degree} of a vertex $i$, denoted $d(i)$, is the cardinality $\card{V(i)}-1$. For example, in the \emph{complete} graph $K_n$ each vertex has degree $n-1$, while in a \emph{cycle} $C_n$ each vertex has degree $2$. The graph $\mathcal{G}$ is said to be \emph{connected} if any two vertices of $\mathcal{G}$ are connected by a path. We denote by $\overline{\mathcal{G}}$ the \emph{complement}
of $\mathcal{G}$, that is, $\overline{\mathcal{G}} = (V, [V]^2\setminus \mathcal{E})$. 
The \emph{disjoint union} of two graphs $\mathcal{G}_1 = (V_1, \mathcal{E}_1)$, $\mathcal{G}_2 = (V_2, \mathcal{E}_2)$, $V_1\cap V_2=\emptyset $, is defined as the graph $\mathcal{G} = (V_1\cup V_2,\mathcal{E}_1 \cup \mathcal{E}_2)$. The \emph{graph join} of $\mathcal{G}_1 = (V_1, \mathcal{E}_1)$ and $\mathcal{G}_2 = (V_2, \mathcal{E}_2)$, denoted $\mathcal{G}_1+\mathcal{G}_2$, is defined as the disjoint union of $\mathcal{G}_1$ and $\mathcal{G}_2$ together with the edges $\{i_1,i_2\}$ for all $i_1\in V_1$ and $i_2\in V_2$. For further background in graph theory, see e.g.  \cite{bondy-murty, diestel, GGL}.

The proofs of the following results are omitted and can be verified by repeated use of Lemma \ref{keylemma}. To illustrate we provide a proof of Proposition \ref{disconnected} which settles the case of disconnected graphs.

\begin{proposition}\label{disconnected}
 Suppose $\mathcal{G}$ is disconnected. Then $\mathcal{G}$ is join-irreducible if and only if 
$\check{\mathcal{G}}$ is isomorphic to the disjoint union of $n$ copies of $K_3$, for some $n\geq 2$.
\end{proposition}
\begin{proof}
 Clearly, if $\check{\mathcal{G}}$ is isomorphic to the disjoint union of $n$ copies of $K_3$, for some $n\geq 2$, then $\mathcal{G}$ is join-irreducible. For the converse, we may suppose  with no loss of generality that $\mathcal G=\check{\mathcal G}$. Let $\mathcal{G}_1$ and $\mathcal{G}_2$ be two connected components of $\mathcal{G}$. Note that $\card{\mathcal{G}_1},\card{\mathcal{G}_2}\geq 2$. Take $i\in \mathcal{G}_1$ and $j\in \mathcal{G}_2$, and let $e=\{i,j\}$. Clearly, $\card{\mathcal{G}_e}=\card{\mathcal{G}}-1$, no vertice is isolated in $\mathcal{G}_e$ and $\mathcal{G}_e$
has one less connected component than $\mathcal{G}$.

Now take $i,i'\in \mathcal{G}_1$ and let $e'=\{i,i'\}$. Clearly, for every such choice of $e'$, we have $e\not \approx e'$. Since $\mathcal{G}$ is join-irreducible, Lemma \ref{keylemma} implies that $\ess f_{\mathcal{G}_{e'}}<\ess f_{\mathcal{G}_e}$. In other words, for every $e'=\{i,i'\}$, $i,i'\in \mathcal{G}_1$, $\ess f_{\mathcal{G}_{e'}}\leq \ess f_{\mathcal{G}}-2$. From Theorem \ref{booleangap} it follows that $\mathcal{G}_1$ must be isomorphic to $K_3$. Since the choice of connected components was arbitrary, we conclude that $\mathcal{G}$ is isomorphic to the disjoint union of $n$ copies of $K_3$, for some $n\geq 2$. 
\end{proof}

To deal with the case of connected graphs, we need to introduce some terminology.
Let $\mathcal{G}=(V, \mathcal{E})$ be a graph. A subset $S\subseteq V$ is said to be \emph{autonomous} if for every $i,i'\in S$ and $j\in V\setminus S$, $i\sim j$ if and only if $i'\sim j$. Moreover, $S$ is said
to be \emph{independent} if for every $i,i'\in S$, $i\not \sim i'$. For simplicity, we refer to autonomous independent sets as $ai$-sets.
We say that $\mathcal{G}$ is \emph{$ai$-prime} if its $ai$-sets are empty or singletons.

\begin{fact} For each $i\in V$, the union of all $ai$-sets containing $i$ is an $ai$-set called the \emph{$ai$-component} of $i$. Moreover, each graph $\mathcal{G}$ decomposes into $ai$-components.
\end{fact}

On the set of $ai$-components of $\mathcal{G}$ there is a graph structure, denoted $\mathcal{G}_{ai}$, in such a way that $\mathcal{G}$ is the lexicographic sum of its $ai$-components and indexed 
by $\mathcal{G}_{ai}$. Note that the graph $\mathcal{G}_{ai}$ is $ai$-prime. 

These constructions are variants of the classical notions of decomposition of graphs and prime graphs (see \cite{ehrenfeucht}).
%Auxiliary lemma A in notes
\begin{lemma}\label{lemma:aux}
 Let $\mathcal{G}=(V, \mathcal{E})$ be a connected graph and suppose that there is $e \in [V]^2\setminus \mathcal{E}$ such that $\mathcal{G}_e$ has no isolated vertices. Then there is $e' \in \mathcal{E}$ such that $\mathcal{G}_{e'}$ has no isolated vertices.
\end{lemma}

Thus, if $\mathcal{G}=(V, \mathcal{E})$ is join-irreducible,  $\mathcal{G}_e$ has an isolated vertex for every $e\in [V]^2\setminus \mathcal{E}$. Moreover,  the nonedge  $e=\{i_1,i_2\}  \in [V]^2\setminus \mathcal{E}$ must be in a $ai$-component or there is $j\in V$ such that $d(j)=2$ and $i_1\sim j\sim i_2$.

We say that a graph $\mathcal{G}=(V, \mathcal{E})$ \emph{satisfies} $(P)$ if for every nonedge  $e=\{i_1,i_2\}  \in [V]^2\setminus \mathcal{E}$ there is $j\in V$ such that $d(j)=2$ and $i_1\sim j\sim i_2$.

Lemma \ref{lemma:aux} and the observation above yield the following.
\begin{corollary}  If  a connected graph $\mathcal G$ is join-irreducible, then $\mathcal G_{ai}$ satisfies $(P)$.
\end{corollary}

Our next proposition describes those graphs satisfying property (P).
\begin{proposition}\label{prop:propP} A  graph $\mathcal{G}=(V, \mathcal{E})$ satisfies $(P)$ if and only if $\mathcal{G}$ is either isomorphic to $K_n$, for some $n\geq 2$, $C_5$, $C_4$ or to a $3$-element path.
\end{proposition} 

\begin{proof}[Sketch proof] We observe that each member of the list satisfies $(P)$. Conversely, suppose that $\mathcal{G}=(V, \mathcal{E})$ satisfies $(P)$. We prove successively:  
\begin{claim}Let $e=\{i_1,i_2\}  \in [V]^2\setminus \mathcal{E}$ and  $j\in V$ such that $d(j)=2$ and $i_1\sim j\sim i_2$. If $e':=\{j, j'\}\in [V]^2\setminus \mathcal{E}$ then either $e_1:=\{i_1,j'\}  \in \mathcal{E}$ and $d(i_1)=2$ or $e_2:=\{i_2,j'\}  \in \mathcal{E}$ and $d(i_2)=2$.
\end{claim}
\begin{claim}\label{claim:clique} Let $i\in V$. If $d(i)\geq 3$, then $\mathcal{G}(i)=(V(i),[V(i)]^2\cap \mathcal{E})$ is isomorphic to $K_n$, for some $n\geq 2$. In fact,
$\mathcal{G}(i)=\mathcal{G}$.
\end{claim}
According to  Claim \ref{claim:clique}, if $\mathcal G$ is not isomorphic to $K_n$, the degree of each vertex is at most $2$. Since $\mathcal G$ satisfies $(P)$, it must be isomorphic to one of the three last members of our list.  
\end{proof}

As a corollary we get the following result. 

\begin{corollary}\label{cor-ai}
If $\mathcal{G}$ is connected and 
join-irreducible, then $\mathcal{G}_{ai}$ is isomorphic to
 $K_n$, for some $n\geq 2$, or  to $C_5$.
\end{corollary}

Clearly,  each $K_n$, $n\geq 2$, and $C_5$ are join-irreducible graphs.
Thus, if $\mathcal{G}=(V, \mathcal{E})$ is an $ai$-prime graph, then $\mathcal{G}$ is join-irreducible if and only if ${\mathcal{G}}$ is isomorphic to $K_n$, for some $n\geq 2$, or to $C_5$.

Now if a connected and join-irreducible graph $\mathcal{G}=(V, \mathcal{E})$ is  not an $ai$-prime graph, then ${\mathcal{G}_{ai}}$ cannot be isomorphic to $C_5$. Indeed, for the sake of contradiction, suppose that ${\mathcal{G}_{ai}}$ is isomorphic to $C_5$. Let $\mathcal{G}_1, \ldots , \mathcal{G}_5$
be the $ai$-components of $\mathcal{G}$ such that $\mathcal{G}_i$ is connected to $\mathcal{G}_{i+1}$, for $i=1,2,3,4$ and $\mathcal{G}_5$ is connected to $\mathcal{G}_{1}$. Assume, without loss of generality, that $\card{\mathcal{G}_1}\geq 2$. Consider $i,i'\in \mathcal{G}_1$, $i_2\in \mathcal{G}_2$ and $i_3\in \mathcal{G}_3$, and let $e=\{i,i_2\}$ and $e'=\{i',i_3\}$. Clearly, $e\not \approx e'$ and $\ess f_{\mathcal{G}_{e}}=\ess f_{\mathcal{G}_{e'}}$. By Lemma \ref{keylemma} it follows that $\mathcal{G}$ is not  join-irreducible which constitutes the desired contradiction. 
Thus, by Corollary \ref{cor-ai} it follows that, in the non $ai$-prime case, if $\mathcal{G}=(V, \mathcal{E})$ is join-irreducible, then ${\mathcal{G}_{ai}}$ is isomorphic to $K_n$, for some $n\geq 2$.

\begin{proposition}Let $\mathcal{G}=(V, \mathcal{E})$ be a connected and non $ai$-prime   graph. Then $\mathcal{G}$ is join-irreducible if and only if ${\mathcal G}$ is isomorphic to one of the following graphs:
\begin{enumerate}[(i)]
\item ${K}_2+\overline{K}_m$,  for some $m\geq 2$;
\item $\overline{K}_n+\overline{K}_m$,
for some $n,m$ with $1\leq n<m$;
\item a graph join $\overline{K}_n+\ldots +\overline{K}_n$ of $r$ copies of $K_n$, for some $r, n\geq 2$.
\end{enumerate}
\end{proposition}
\begin{proof}[Sketch proof] As observed if $\mathcal{G}=(V, \mathcal{E})$ is join-irreducible, then ${\mathcal{G}_{ai}}$ is isomorphic to $K_r$, for some $r\geq 2$.  For $r=2$, it is clear that  $\mathcal G$ is isomorphic to $\overline{K}_n+\overline{K}_m$ for some $n, m\geq 1$. For $r=3$, if there is $i\in V$ such that $d(i)=2$, then we show that $\mathcal G$ is isomorphic to ${K}_2+\overline{K}_m$,  for some $m\geq 2$. If $r \geq 3$ and for every $i\in V$ we have $d(i)>2$, then we show that  $\mathcal G$ is the join $\overline{K}_n+\ldots +\overline{K}_n$  of $r$ copies of ${K}_n$ for some $n\geq 2$.
\end{proof}

From these results, we obtain the description of the join-irreducible graphs.

\begin{theorem} \label{thm:jigraphs}Let $\mathcal{G}=(V, \mathcal{E})$ be a graph. Then $\mathcal{G}$ is irreducible if and only if $\check{\mathcal G}$ is isomorphic to one of the following graphs:
\begin{enumerate}[(i)]
 \item a disjoint union of $n$ copies of $K_3$, for some $n\geq 2$;
\item $C_5$;
\item ${K}_2+\overline{K}_m$,  for some $m\geq 2$;\item $K_n$, for some $n\geq 2$;
\item $\overline{K}_n+\overline{K}_m$,
for some $n,m$ with $1\leq n<m$;
\item a graph join $\overline{K}_n+\ldots +\overline{K}_n$ of $r$ copies of $K_n$, for some $r, n\geq 2$.
\end{enumerate}
\end{theorem}

\end{document}